\documentclass[12pt]{article}

\usepackage{fullpage,amsfonts,amsmath, amssymb,mathrsfs,graphicx}

\usepackage{verbatim,url}

\usepackage[top=2.54cm, bottom=2.54cm, left=2.54 cm, right=2.54 cm]{geometry}
\newcommand{\lab}[1]{\label{#1}}                

\usepackage[usenames,dvipsnames]{color}

\newcommand{\jcoma}[1]{{\bf [~jane 07/03/17:\ } {\color{blue} #1}{\bf~]}}
\newcommand{\jcomb}[1]{{\bf [~jane 20/03/17:\ } {\color{Cyan} #1}{\bf~]}}
\newcommand{\jtb}[1]{{\color{Purple} #1}}

\newcommand{\jcome}[1]{{\bf [~jane 19/04/17:\ } {\color{Cyan} #1}{\bf~]}}

\newcommand{\nac}[1]{{\color{cyan}{\bf [Nick:}} {\color{ForestGreen}{\em #1}}{\color{cyan}{\bf ]}}}

 \newcommand{\nek}[1]{{\color{red} #1}}

\newcommand{\remove}[1]{}
\newcommand\eqn[1]{(\ref{#1})}

\newcommand{\be}{\begin{equation}}
\newcommand{\bel}[1]{\begin{equation}\lab{#1}\ }
\newcommand{\ee}{\end{equation}}
\newcommand{\bea}{\begin{align}}
\newcommand{\eea}{\end{align}}
\newcommand{\bean}{\begin{align*}}
\newcommand{\eean}{\end{align*}}
\newcommand{\ceil}[1]{\lceil #1 \rceil}
\newcommand{\floor}[1]{\lfloor #1 \rfloor}

\newtheorem{thm}{Theorem}

\newtheorem{con}[thm]{Conjecture}
\newtheorem{lemma}[thm]{Lemma}

\def\proof{\noindent{\bf Proof.\ }  }
\def\qed{~~\vrule height8pt width4pt depth0pt}

\def\M{{\mathcal M}}

\def\ex{{\mathbb E}}
\def\pr{{\mathbb P}}

\def\eps{\epsilon}

\def\ss{\smallskip}

\def\no{\noindent}

\renewcommand{\ge}{\geqslant}
\renewcommand{\le}{\leqslant}
\renewcommand\epsilon{\varepsilon}
\renewcommand\emptyset{\varnothing}

\def\eref#1{$(\ref{#1})$}
\def\sref#1{\S$\ref{#1}$}
\def\lref#1{Lemma~$\ref{#1}$}
\def\tref#1{Theorem~$\ref{#1}$}
\def\cjref#1{Conjecture~$\ref{#1}$}

\def\fref#1{Figure~$\ref{#1}$}

\date{}

\title{
Full rainbow matchings in graphs and hypergraphs
}

\author{Pu Gao\thanks{Research supported by the ARC grant DE170100716.} \\
{\small School of Mathematics}\\[-0.8ex]
{\small Monash University}\\
{\small \texttt{jane.gao@monash.edu}} \and 
Reshma Ramadurai\\
{\small School of Mathematics \& Statistics}\\[-0.8ex]
{\small Victoria University of Wellington}\\
{\small \texttt{Reshma.ramadurai@vuw.ac.nz}} \and 
Ian M.\ Wanless\thanks{Research supported by the ARC grant DP150100506.}\\
{\small School of Mathematics}\\[-0.8ex]
{\small Monash University}\\
{\small \texttt{ian.wanless@monash.edu}} \and 
Nick Wormald \thanks{Research supported by the Australian Laureate Fellowships grant FL120100125.}\\
{\small School of Mathematics}\\[-0.8ex]
{\small Monash University}\\
{\small \texttt{nick.wormald@monash.edu}}
}

\begin{document}
\maketitle

\begin{abstract}

Let $G$ be a simple graph that is properly edge coloured with $m$ colours and let $\M=\{M_1,\ldots, M_m\}$ be the set of $m$ matchings induced by the colours in $G$. Suppose that $m\le n-n^{c}$, where $c>9/10$, and every matching in $\M$ has size $n$. Then $G$ contains a full rainbow matching, i.e.\ a matching that contains exactly one edge from $M_i$ for each $1\le i\le m$. This answers an open problem of Pokrovskiy and gives an affirmative answer to a generalisation of a special case of a conjecture of Aharoni and Berger.

Related results are also found for multigraphs with edges of bounded
multiplicity, and for hypergraphs.

Finally, we provide counterexamples to several conjectures on full rainbow matchings made by Aharoni and Berger.
\end{abstract}

\section{Introduction}

Throughout this paper the setting is a multigraph $G$ whose edges are properly coloured with $m$ colours, so that each colour $i$ induces a matching $M_i$. We say that $G$ contains a {\em full rainbow matching} if there is a matching $M$ that contains exactly one edge from $M_i$ for each $1\le i\le m$. This paper is motivated by the following conjecture of Aharoni and Berger \cite[Conj.~2.4]{ahbe}.

\begin{con}\label{cj:AB}
If $G$ is bipartite and each matching $M_i$ has size $m+1$ then $G$ has a full rainbow matching.
\end{con}

Consider a $k\times n$ array $A$. A {\em partial transversal of length
  $\ell$} in $A$ is a selection of $\ell$ cells of $A$ from different
rows and columns, and containing different symbols.  A {\em
  transversal} of $A$ is a partial transversal of length $\min(k,n)$.
If $A$ has no repeated symbol within a row it is called {\em
  row-Latin}.  We say that $A$ is {\em Latin} if it and its transpose
are both row-Latin. If $k=n$ and $A$ is Latin and contains exactly $n$
symbols, then $A$ is a {\em Latin square}.  \cjref{cj:AB} was
motivated by a longstanding conjecture of Stein \cite{Ste75}, that
every $(n-1)\times n$ row-Latin array has a transversal. This is
equivalent to the restriction of \cjref{cj:AB} to the case where each
matching covers the same set of vertices on one side of the bipartite
graph. Stein's conjecture was in turn motivated by the question of
what length of partial transversal can be guaranteed to exist in a
Latin square. His conjecture implies that every Latin square of order
$n$ has a partial transversal of length $n-1$ (this statement was
independently conjectured by Brualdi slightly earlier; see
\cite{transurv} for a full survey of these conjectures and related
results). The best result to date is by Hatami and Shor \cite{HS08},
who showed that every Latin square of order $n$ has a partial
transversal of length $n-O(\log^2 n)$.  It is known that for even
orders $n$ there are at least $n^{n^{3/2}(1/2-o(1))}$ (equivalence
classes of) Latin squares that do not have transversals
\cite{CW17}. However, a famous conjecture of Ryser \cite{ry} states
that all Latin squares of odd order have transversals. In terminology
similar to \cjref{cj:AB}, Ryser's conjecture is that if $G$ is
$K_{m,m}$ and $m$ is odd, then $G$ should have a full rainbow
matching.  This conjecture is known to fail if a single edge is
removed from $K_{m,m}$. Also there are Latin arrays of odd order $n$
containing more than $n$ symbols but having no transversal (again, see
\cite{transurv} for details).

Bar{\'a}t and Wanless \cite{BW14} considered an intermediate step
between \cjref{cj:AB} and its variant (that we know fails) with $m+1$
replaced by $m$. They showed that $\lfloor m/2\rfloor-1$ matchings
of size $m+1$ together with $m-\lfloor m/2\rfloor+1$ matchings
of size $m$ need not have a full rainbow matching. They also constructed
$m$ matchings of size $m$ inducing a bipartite multigraph with $m$ vertices 
in the first part of the bipartition and $m^2/2-O(m)$
in the second part, and with no rainbow matching. This raises the question
of how large one part can be before a rainbow matching is unavoidable.
In \cite{BHWWW17} it is shown that if one part has $m$ vertices and the other
has at least $\big\lceil\frac{1}{4}(5-\sqrt{5})m^2\big\rceil$ vertices 
then there will be a full rainbow matching. 
Clearly, the threshold is quadratic in $m$ for this problem. However,
things change significantly if the induced bipartite graph must be simple.
Montgomery, Pokrovskiy and Sudakov \cite{MPS18}, 
showed in that case that if one part has $m$ vertices and the other has 
at least $\epsilon m^2$ vertices then there will be many full rainbow
matchings.  Also Keevash and Yepremyan \cite{KY18} showed that if one
part has $m$ vertices and the other has at least $m^{399/200}$
vertices then there will be a full rainbow matching. In particular,
the threshold for this variant of the problem is subquadratic.

The first progress towards \cjref{cj:AB} was by Aharoni, Charbit and
Howard \cite{ahchho}, who showed that $n$ matchings of size
$\lfloor7n/4\rfloor$ must have a full rainbow matching. 
The $\lfloor7n/4\rfloor$ term was successively improved to
$\lfloor5n/3\rfloor$ by Kotlar and Ziv \cite{kozi},
then $(3+\epsilon)n/2$ by Clemens and Ehrenm\"{u}ller \cite{cleh}, and
then $\lceil3n/2\rceil+1$ by Aharoni, Kotlar and Ziv \cite{ahkozi}.
Finally, for any fixed $\eps>0$, Pokrovskiy \cite{Pok17} showed that if
the matchings are edge-disjoint (so that $G$ is simple) then
$n$ matchings of size $(1+\epsilon)n$ have a full rainbow matching when $n$
is sufficiently large. He also posed
two challenges regarding improving the error term in his result (we
make some progress in this direction) and generalising to bipartite
multigraphs (a feat he himself achieved in \cite{Pok16}).

A related result is due to H\"aggkvist and Johansson \cite{hajo},
who showed that if the matchings are all perfect and edge-disjoint
then $n$ matchings of size $(1+\epsilon)n$ can be decomposed into 
full rainbow matchings, provided $n$ is sufficiently large. 

All the results discussed so far pertain to bipartite graphs. So far, this case has attracted more scrutiny than the
unrestricted case. However, Aharoni {\em et al.}~\cite{ABCHS} and Bar\'at,
Gy\'arf\'as and S\'ark\"ozy~\cite{BGS}
both consider the question of how large a
rainbow matching can be found across any set of matchings.  
The former paper makes a conjecture which includes this variant of
\cjref{cj:AB} as a special case:
\begin{con}\label{cj:ABCHS}
If each matching has size $m+2$ then $G$ has a full rainbow matching.
\end{con}
It is not viable to replace $m+2$ by $m+1$. For example, a
$1$-factorisation of two copies of $K_4$ provides $3$ matchings of
size $4$ that do not possess a full rainbow matching.

Our aim is to investigate approximate versions of \cjref{cj:ABCHS}.
Our main result is an analogue of Pokrovskiy's Theorem from
\cite{Pok17}, but without the requirement that $G$ is bipartite. Our
results in this direction are stated in the next section, and include
some information on the role of the maximum degree of the graph
$G$. Related to this, in \sref{s:CEs}, we discuss and refute several
further conjectures on full rainbow matchings made by Aharoni and
Berger \cite{ahbe} and Aharoni {\em et al.} \cite{ABCHS}. In
particular, both papers include the following conjecture, which turns
out not to be true.

\begin{con}\label{cj:falseconjbip}
  Let $G$ be a bipartite multigraph, with maximum degree $\Delta(G)$,
  whose edges are (not-necessarily properly)
  coloured. If every colour appears on at least $\Delta(G)+1$ edges, then
  $G$ has a full rainbow matching. 
\end{con}



Finally, we note that the main result in a recent preprint of Keevash and Yepremyan
\cite{KY17} implies  that in any multigraph with edge multiplicities
$o(n)$ that is properly edge-coloured by $n$ colours with at least 
$n(1+\epsilon)$ edges of each colour, there must be a rainbow matching
that is close to full.

\section{Main results}


Recall our setting:  $G$ is a multigraph that is properly edge coloured with $m$ colours, and matching $M_i$ is induced by colour $i$. Let $\M=\{M_1,\ldots, M_m\}$ denote the family of $m$ matchings. We say $\M$ is {\em non-intersecting} if $G$ is a simple graph, i.e. if $M_i\cap M_j =\emptyset$ for every $1\le i<j\le m$. Otherwise, it is called intersecting. Let $\Delta(\M)=\Delta(G)$ denote the maximum degree of $G$.
 
\begin{thm}\lab{thm:simple} 
Suppose that $0\le \delta<1/4$ and $0<c<(1-4\delta)/10$ and $n$ is sufficiently large. If $\M$ is a non-intersecting family of $m\le (1-n^{-c})n^{1+\delta}$ matchings, each of size $n$, and $\Delta(\M)\le (1-n^{-c})n$, then $\M$ contains a full rainbow matching. 
\end{thm}



Clearly, for any $\M$ we always have $\Delta(\M)\le m$. Thus by taking $\delta=0$, we immediately have the following corollary.

\begin{thm}\lab{thm2:simple} 
Suppose that $0<c<1/10$ and $n$ is sufficiently large. If $\M$ is a non-intersecting family of $m\le (1-n^{-c})n$ matchings, each of size $n$, then $\M$ contains a full rainbow matching.
\end{thm}

{
\no{\bf Remark:} \tref{thm2:simple} proves an approximate version of
\cjref{cj:ABCHS}, and thus of \cjref{cj:AB}. Compared
with~\cite{Pok16} we get better approximation by improving $m$ from
$(1-o(1))n$ to $n-n^{9/10+\eps}$ and of course \tref{thm2:simple} also
approximates the non-bipartite case of \cjref{cj:ABCHS}. However
\tref{thm2:simple} does not cover the multigraph case,
whereas~\cite{Pok16} does.  \ss
}



\no {\bf Remark:} 
One can ask about possible strengthenings of this theorem. For
instance, Aharoni and Berger~\cite[Conjecture~2.5]{ahbe} conjectured
essentially that one can drop the upper bound on $m$, and even drop
the condition that each element of $\M$ is a matching, as long
as $\Delta(\M)\le n-1$. It turns out that their conjecture
is false, as shown by a graph whose components are double stars, which
we give in \sref{s:CEs}. But the question for matchings remains
open.


Our proof for \tref{thm:simple} easily extends to intersecting $\M$ where the underlying graph $G$ is a multigraph with relatively low multiplicity. 

\begin{thm}\lab{thm:multiple}
  For every $\eps_0>0$, if $\M$ is a family of $m\le (1-\eps_0) n$ matchings, each of size $n$, and every edge is contained in at most $\sqrt{n}/\log^2 n$ matchings, then $\M$ contains a full rainbow matching. 
\end{thm}

The proof of \tref{thm:simple} also immediately extends to rainbow matchings in uniform hypergraphs. A hypergraph $G=(V,E)$ is defined on a set of vertices $V$ where the set $E$ of hyperedges is a set of subsets of $V$. We say $G$ is a $k$-uniform hypergraph if every hyperedge has size $k$.  A matching $M$ in $G$ is a set of hyperedges such that no vertex in $G$ is contained in more than one hyperedge in $M$.

\begin{thm}\lab{thm:hypergraph}
  For every $\eps_0>0$ and every integer $k\ge 2$, if $\M$ is a family of $m\le (1-\eps_0) n$ edge disjoint matchings in a $k$-uniform hypergraph $H$, each of size $n$, and every pair of vertices is contained in at most $\sqrt{n}/\log^2 n$ hyperedges, then $\M$ contains a full rainbow matching. 
\end{thm}
\no {\bf Remark:} If we restrict further by saying that every pair of vertices is contained in at most a constant number of hyperedges, then \tref{thm:hypergraph} holds if we replace $m\le (1-\eps_0)n$ by $m\le (1-n^{-1/10+\eps_0})n$. That is, we can fully recover \tref{thm2:simple} for arbitrary $k\ge 2$. This indeed covers some interesting families of hypergraphs such as linear hypergraphs where no two vertices are contained in more than one hyperedge. By further restricting the maximum degree of $G$, improved bounds on $m$ can be achieved with minor modification of the proof. Similarly the bound on $m$ in \tref{thm:multiple} can   easily be improved by restricting to smaller multiplicity, or placing additional constraints on the maximum degree. We will not develop this idea further in this paper.


 \section{A heuristic approach} 
 \label{sec:heuristic}

 In this section, we give a simplified description of the algorithm we use for full rainbow matchings, and also a heuristic argument as to why we expect it to successfully find the full rainbow   matching  required for  \tref{thm:simple}. The actual proof, showing that all aspects work as intended,  will be given in \sref{sec:proof}.
 
Given a nonintersecting family $\M$ of $m$ matchings of size $n$, we do the following.

First, randomly partition $\M$ into subfamilies containing about $\eps
m$ matchings each, which we call ``chunks". Here $\eps$ is a function
of $n$.

Next, ``process'' the chunks iteratively using the following three
steps in iteration $i$.

\begin{itemize}

\item[(i)] Pick one edge u.a.r.\ from each of the matchings in chunk
  $i$. Any picked edge $x$ that is not incident to any other is added
  to the rainbow matching $M_0$ that will be outputted, and the end
  vertices of $x$ are deleted from the graph. Edges that ``collide''
  with others are not added (but see step (iii)).

\item[(ii)] For each vertex surviving the first step, calculate the
  probability that it was deleted, and then artificially delete with
  such a probability to ensure that all vertices have the same
  probability of surviving these first two steps of the iteration. (A
  suitable probability will be specified in the precise analysis of
  the algorithm.)

\item[(iii)] For any matching $M$ containing a ``colliding'' edge in
  step (i), greedily choose a replacement edge $x$ in $M$ to add to $M_0$
  and delete the end vertices of $x$ from the graph.

\end{itemize}
 
\noindent These steps are performed for all chunks except the
last, which is treated instead by greedily choosing edges from the
remaining matchings. (We will show that this is highly likely to
succeed.)

We now give a rough overview of the analysis of the algorithm.
Let $\tau$ denote the number of iterations of the algorithm, i.e.\ the number of chunks. (For definiteness, we call the treatment of the last chunk an
``iteration'', even though it is treated differently to the other chunks.)
Also, let $d_v^{j} (i)$ denote the number of edges of chunk $j$ that are (still) incident with vertex $v$ after iteration $i$. We will specify  functions $g$ and $r$, and the correct probabilities in  step (ii), for the following to hold iteratively for each $i=1,\ldots,  \tau-1 $ and all $j > i$:
\begin{itemize}
\item[(a)]   
  after iteration $i$, all  surviving  matchings have size approximately $r(i\eps)n$;
  
\item[(b)]  $d_v^{j} (i) \approx \eps g(i\eps) d_v$  for all surviving vertices $v$, where 
  $d_v$ is the degree of $v$ initially (in $G$).
  
\end{itemize}
 
\noindent
Here the sign $\approx$ is used to  denote some version of  concentration around the stated  value.

Specific versions of (a) and (b) are proved using the initial
concentration, together with an inductive argument that computes the
expected changes in the variables concerned during each iteration and
shows concentration close to the expected changes. The precise
inductive statements are chosen with a margin of error that rigorously
contains the effect of collisions, so we will ignore these in the
outline here.  Our estimates are least accurate near the end of the
algorithm, which is why we use a greedy algorithm at that stage. This
final iteration works because there are few matchings left compared to
their remaining sizes.

Here is an outline of why we expect the algorithm to succeed. Initially we have $|M|= n$ for all matchings $M$.  Since the matchings are randomly allocated into chunks, we expect the initial setup to satisfy
$$d_v^{j}(0) \approx \eps d_v, \quad \mbox{for all $j\ge 1$ and $v\in G$.}$$
Assume that the first $i$ iterations of the algorithm are complete. 
 For iteration $i+1$,  we specify the probability  of $v$ being artificially deleted in step (ii) so that the probability of surviving the first two steps is $f(i\eps)$ for every remaining vertex.  Then for every vertex $v$ that survives iteration $i+1$, the expected change in $ d_v^{j} (i)$ ($j>i$) is
\be\lab{drec}
\ex ( d_v^{j} (i+1) -d_v^{j} (i)) \approx  -f(i\eps)  d_v^{j} (i),
\ee
since each of $v$'s neighbours is deleted with probability $f(i\eps)$.  (Here we ignore collisions and the effect of step (iii), as mentioned above.)
Hence, the degrees of two different vertices remain roughly in the
same proportion as long as they both survive.

For each matching in chunk $j$, the expected change in its size while chunk $i$ is processed is roughly  
\bel{Mrec}
-2f(i\eps)|M(i)|,
\ee 
neglecting what will turn out to be an $O(\eps^2n)$ error from the case that the two ends of an edge in $M$ are both deleted.

Consider a vertex $v$ that has survived $i$ iterations. During iteration $i+1$, each matching in  chunk $i+1$  has one of its edges chosen. So, using (a) above, any given edge 
in a matching in chunk $i+1$ that is incident with $v$
is chosen with probability $p\approx 1/\big(r(i\eps)n\big)$. 
This means that the probability that $v$  is not deleted in step (i) of iteration $i+1$ is roughly $(1-p)^{d_v^{{i+1}} (i) }\approx 1-pd_v^{{i+1}} (i) $, and hence the probability that it is deleted here is roughly 
$$
\frac{ d_v^{i+1} (i)}{  r(i\eps)n}.  
$$ 
With this in mind, we can define  $f(i\eps)$ so that it is approximately equal to the maximum value of this probability over all $v$, which is determined by the maximum vertex degree. (We also add a little elbow-room to account for the collisions.) Then the probability 
can be appropriately specified in step (ii). We can find the maximum value by tracking the maximum degree via~\eqn{drec}.  Let $\gamma  n$ denote the maximum degree of the initial graph $G$.   By our assumption in Theorem~\ref{thm:simple}, $\gamma\le 1-n^{-c}$. Using (b), we find that the approximate size of $\max_{v} d_v^j(i)$ is $\eps g(x)\gamma n$,
    where $x=i\eps$. We have 
$$
f(x) \approx \frac{\eps{  g(x)\gamma} n}{  r(x){n} } =\frac{\eps\gamma g(x)}{r(x)}.
$$

Letting $\hat f(x) =\gamma  \frac{ { g}(x)}{  r(x) }$, we have $f(x)\approx\eps\hat f(x)$. Then, if the size of each matching is approximated, as mentioned in (a), by $r(i\eps)n$, equation~\eqn{drec}    suggests (as $\eps\to 0$, and applying it to a vertex of maximum degree)  the differential equation  
$$
{g}'(x)=-\hat f(x){g}(x) = -\gamma\frac{  { g}(x)^2}{  r(x) }.
$$
Similarly,~\eqn{Mrec} suggests
$$
r'(x) = -  2\gamma g(x).
$$
But then $dg/dr = g/ 2r$ which gives $r=Cg^2$. Initially, $r(0)=1$ and $g(0)=1$ which yields $C=1$. So, the solution to these differential equations is
\bel{desol}
r(x) =\left (1-\gamma x \right)^2, \quad g(x)= 1-\gamma x.
\ee
Thus, we have $d_v^j(i)\approx \eps (1-\gamma i\eps) d_v$ for every surviving vertex $v$. For those vertices whose degrees are initially lower than the maximum degree, the derivative is proportionally lower, and hence the degrees stay in proportion. The process cannot `get stuck' until the error in the approximation in assumption (a) becomes significantly large compared  to $r(x)$. The function $r(x)$ is positive for all $0\le x\le 1$, because our hypothesis $\Delta({\mathcal M})\le (1-n^{-c})n$ guarantees that $1-\gamma x\ge n^{-c}>0$.  
Thus, if the error of approximations is small enough (as we shall show),  the process proceeds until the last iteration. The final (greedy) iteration is shown to work using the previous analysis to estimate the size of the remaining matchings.  Note that $g(x)$ tends to zero along with $r(x)$.
 
%
 
The astute reader may have noticed that, in this above sketch of proof,  all we needed  from assumption (b) was an upper bound on all vertex degrees, and this is the approach we will take in the formal proof in the following section. Thus, we will replace assumption (b) by 

(b$'$) $d_v^j(v)\le \eps {g(i\eps) \gamma n}$ for all $j> i$ and all surviving vertices $v$.

\section{Algorithm and proof}
~\lab{sec:proof}

In this section, we define  the algorithm precisely  and then analyse it to prove Theorem~\ref{thm:simple}.

\subsection{The algorithm}

Let $\M = \{M_1, M_2, \ldots, M_m\}$ be a non-intersecting family of matchings, each of size $n$.  
 The algorithm has an initial stage, then some repeated iterations, then one final iteration.
 The initial stage consists of the following.   First order the matchings in  $\M$ uniformly at random (u.a.r.). Then, for some $\eps>0$ of our choosing,   partition $\M$ into ``chunks"  $\M^1, \M^2, \ldots$ where $\M^1$ contains the first $\lceil \eps {m}\rceil$ matchings, $\M^2$ contains the next $\lceil \eps {m}\rceil$  matchings, and so on, up to the last chunk, which   contains at most $\lceil\eps m\rceil$ matchings. For ease of calculations, we will choose $\eps$ so that $\eps {m}$ is an integer.

We next define some notation useful in defining the iterations of the algorithm. Let $G$ be the graph induced by $\cup_{M \in \M}M $, and let $V$ be its vertex set. For any $u\in V$,   let $E_u$ denote the set of edges in $G$ that are incident with $u$. 
During the algorithm, vertices are removed from consideration for several distinct reasons, which we discuss shortly.
The set $U(i)$ is the set of vertices that were removed at some point during the first $i$ iterations. After $i$ iterations, vertices in $V\setminus U(i)$ are said to be \emph{surviving} and matchings are said to be \emph{surviving} if they do not belong to the first $i$ chunks.
Edges are said to be \emph{surviving} if both their endpoints are surviving vertices and they are part of a surviving matching.
At any point in the algorithm $M_0$ denotes the set of edges added so far to the rainbow matching (initially $M_0=\emptyset$).  The graph $G(i)$ denotes the graph with vertex set $V\setminus U(i)$ and edge set $\bigcup_{j>i} \M^j$ restricted to $V\setminus U (i) $. For all matchings $M \in \mathcal{M}$, we let $M(i)$ denote $M \cap E(G (i))$ and $\M^{j}(i)= \{M(i): M\in \M^j\}$.  Let $E(\mathcal{M}^j(i))=\cup_{M\in \M^j} M{(i)}$, so that $E(\mathcal{M}^j(i))$ is the set of edges in matchings in chunk $j$ that still survive after iteration $i$. The above definitions are all intended to apply to the $i=0$ case in the obvious way, with $U(0)=\emptyset$, $G(0)=G$, and so on.

After the initial stage, 
the algorithm performs iterations consisting of the three steps below. 
We consider the situation after $i\ge0$ iterations have been completed,
and describe how to perform the $(i+1)$-st iteration.
 For simplicity, we describe certain edges and vertices being deleted from $G(i)$ as the algorithm progresses. More accurately,  
  the algorithm takes a copy of $G(i)$ at the start of the $(i+1)$-st iteration and edits this copy, which will end up becoming $G(i+1)$. 


We assume that $f$ is a given function (and will specify a
particular one below).

\begin{itemize}
\item[(i)]
  For each $M \in \M^{i+1}(i)$, choose one edge in $M$ u.a.r.. Let $\Psi(i+1)$ denote the set of edges that are chosen. Vertices incident with edges in $\Psi(i+1)$ are called \emph{marked}.
For   ${x \in \Psi(i+1)} $, if $x\cap {y}\neq \emptyset$ for some  ${y\in \Psi(i+1)}$, we say there is a {\em vertex collision} involving $x$. For each  ${x \in \Psi(i+1)} $ not involved in such a collision, add $x$ into $M_0$ and delete the end vertices of $x$ from $G(i)$; vertices deleted this way are called ``killed''.
 
 \item[(ii)] Independently delete each existing vertex $v$ in $G(i)$ with probability $P_{i+1}(v)$ where
 \[
 Q_{i+1}(v)+P_{i+1}(v)(1-Q_{i+1}(v)) = f(i\eps),
 \]
 and $Q_{i+1}(v)$ denotes the probability that $v$ is marked in
 step (i).  Vertices deleted this way are called ``zapped''.  If
 $P_{i+1}(v)<0$ or $P_{i+1}(v)>1$ for some $v$ then restart the algorithm.
 
 \item[(iii)]   Deal with vertex collisions greedily. Let $\Phi(i+1)$ denote the set of matchings in chunk ${i+1}$ that are not processed yet due to a vertex collision in step (i). Sequentially for each $M\in\Phi(i+1)$, choose a valid edge $x\in M$ using a greedy algorithm; e.g.\ choose $x$ incident with a vertex with the lowest index. Add $x$ into the rainbow matching $M_0$ and delete the end vertices of $x$ from the remaining graph. Unmark any vertices that were marked but not deleted. 
 \end{itemize}
 
The final iteration of the algorithm consists of treating the last
  chunk of matchings. Here edges are chosen greedily one by one from
  those matchings. A simple observation is as follows. If we choose
an edge $x$ that can validly be added to $M_0$, then the removal of the end
vertices of $x$ will decrease the size of each remaining matching by
at most 2. Hence, when the algorithm comes to process the last
  chunk, if the sizes of the remaining matchings are all at least
twice the number of matchings remaining, then a full rainbow
  matching will be successfully completed by the greedy method.

We repeat the following definitions from
  \sref{sec:heuristic}. Let $\gamma n$ denote the maximum degree of
  the initial graph $G$, let $d_v^{j} (i)$ denote the number of edges
  of chunk $j$ that are (still) incident with vertex $v$ after
  iteration $i$, and let $\tau$ denote the number of iterations of the
  algorithm. We have $\tau = \lceil m/( \epsilon {m})\rceil=
  \lceil 1/\eps\rceil$.

\subsection{Proof of {\boldmath\tref{thm:simple}}}
\label{s:proof}
We first change the definition of $\gamma$ slightly from \sref{sec:heuristic}: from now on, set $\gamma=1-n^{-c}$.  
By the hypotheses of the theorem, we may assume that $\M$ is non-intersecting and contains $m=\lfloor\gamma n^{1+\delta}\rfloor$  matchings each of size $n$, and $\Delta(\M)\le \gamma n$. 
(This exact value of $m$ is achieved by adding, if necessary, new matchings that are vertex-disjoint from all previous ones. This does not affect $\Delta(\M)$ or the existence of a full rainbow matching.)  Let $\eps>0$ be a function of $n$, to be specified later,  such that $n^{-1/3}<\eps=o(1)$. Recall that $\tau= \lceil1/\eps  \rceil$ is the number of iterations of the algorithm. It must be noted that our randomised algorithm only applies to the first $\tau-1$ iterations.

For simplicity, we let $r_i$ and $g_i$ denote $r(i\eps)$ and $g(i\eps)$, respectively, where $r(\cdot)$ and $g(\cdot)$ are given in \eref{desol}.
For $0\le i\le {\tau-1}$, we will specify non-negative real numbers $a_i$ and $b_i$ such that at the start of the $(i+1)$-st iteration of the algorithm, the following hold with probability $1-o(\eps)$:  
 \begin{enumerate}
\item[(A1)] every surviving matching has size between $ r_in-a_i$ and $ r_in+a_i$, and
\item[(A2)] every surviving vertex $v$ satisfies
\[
d_v^j(i) \le  \eps {\gamma g_i n }+b_i, \quad \mbox{for all $j > i$}.
\]
\end{enumerate} 
Values of the function $f$ required in step (ii) of the $(i+1)$-st 
iteration will be defined   by    
 \bel{feq}
 f(i\eps) ={\eps\gamma\frac{g_i}{r_i}}+c_i 
 \ee
 where $c_i\ge 0$ will also be specified. 
\remove{
For sequences $f_n$ and $g_n$, we say $f_n=O(g_n)$ if $|f_n|\le C|g_n|$ and all $n$, and for some absolute constant $C$ that is independent of $n$, $\eps$ and $i$. Recall that $\eps$ is a function of $n$ to be specified later. For now we only assume that 
\be
\eps\ge n^{-1/3}.
\ee
}



The proof is by induction. For the base case, $i=0$, we regard the state at the start of the first iteration. The initial graph is $G=G(0)$ and we define  $a_0=0$ and $b_0=(\eps {\gamma n})^{ 1/2}\log {n}$. Then (A1) is trivially true. To verify  (A2) holds for $i=0$, we need to consider the variation in degrees   caused by the initial random permutation of the matchings.

\begin{lemma}\lab{lem:base}
With probability $1-o(\eps)$, property {\rm (A2)} holds for $i=0$. 
\end{lemma}
 
\newcommand{\med}{\lambda}

\proof Here $d_v^j(0)$ is determined by the random permutation $\pi$ of matchings in $\M$.  Obviously $\ex d_v^j(0)=\eps d_v \le \eps \gamma {n}$.
We will apply McDiarmid's inequality~\cite[Theorem 1.1]{McDiarmid}
to prove concentration. Let $\med$ denote the median of $d_v^j(0)$. Observe that 
\begin{itemize}
\item interchanging two elements in $\pi$ can affect $d_v^j(0)$ by at most $\varrho=1$, because all edges incident with $v$ in $G(0)$ belong to different matchings;
\item for every $s>0$, if $d_v^j(0)\ge s$ then there is a set of $s$ elements $\{i_1,\ldots, i_{s}\}\subseteq [m]$ such that $\pi(i_1),\ldots, \pi(i_{s})$ certifies $d_v^j(0)\ge s$. 
\end{itemize}  
By McDiarmid's inequality, for any $t\ge 0$,
\be
\pr(|d_v^j(0)-\med|\ge t) \le 4 \exp\left(-\frac{t^2}{16(\med+t)}\right).\label{eq:McDiarmid}
\ee
It follows immediately that
\begin{align*}
|\ex d_v^j(0) - \lambda| &\le \ex |d_v^j(0)-\lambda| \le \int_{t=0}^{\infty} 4 \exp\left(-\frac{t^2}{16(\med+t)}\right) dt\\
& \le \int_{t=0}^{\med} 4 e^{- t^2/32 \med } dt +\int_{t=\med}^\infty 4  e^{ -t  / 32 } dt=O(\sqrt{\lambda}+1).
\end{align*}
This implies that $\med = \ex d_v^j(0) +O\Big(1+\sqrt{\ex d_v^j(0)}\Big)$. Since $\ex d_v^j(0)\le \eps\gamma n\to\infty$ as $n\to\infty$, we have $\med\le \eps\gamma n+O(\sqrt{\eps\gamma n})$.  Hence,~\eqn{eq:McDiarmid} with $t=\frac12 \sqrt{\eps\gamma n }\log {n}=b_0/2$  yields
\begin{align*}
\pr\big(|d_v^j(0) - \med| \ge b_0/2) & =\exp(-\Omega(\log^2 n)),
\end{align*}
since $16(\med+t )=O( \eps \gamma n)=O(t^2/\log^2 n)$.  As $ \med\le  \eps\gamma n+o(b_0)$, this means  $\pr( d_v^j(0)   \ge  \eps\gamma n + b_0)  =o(n^{-6})$.
Taking union bound over the $O(mn)=O(n^{2+\delta})$ choices for $v$ and $O(m)=O(n^{1+\delta})$ choices for $j$, we can conclude that with probability $1-o(\eps)$ we have $d_v^j(0)\le \eps\gamma n + b_0$ for every $v$ and $j$. 

 This verifies (A2) for $i=0$. \qed\ss

Next assume the claim holds for some $i\ge 0$, i.e., we assume (A1) and (A2) hold after the first $i$ iterations of the algorithm.  Note that most edges in $\Psi(i+1)$ will have their endpoints killed in step (i) whereas some will survive due to vertex collision. 
We say a vertex is  {\em condemned} if it is either zapped or marked. We desire each vertex to be condemned with probability $f(i\eps)$ 
as specified in~\eqn{feq}.
This is made use of in step (ii). Of course, we require that $0\le P_{i+1}(v)\le 1$, which is true if $Q_{i+1}(v) \le f(i\eps)\le 1$. By (A1), after the $i$-th iteration  every surviving matching has size at least $r_in -a_i$, which implies that the probability of a given edge being chosen is at most $1/(r_in-a_i)$.  From (A2), the degree of a vertex is at most $\eps {\gamma g_i n}+b_i$, so we have
$$Q_{i+1}(v) \le \frac{\eps {\gamma g_i n} +b_i}{r_i n-a_i}.$$
Hence, $Q_{i+1}(v)\le f(i\eps)$ would be guaranteed by
\bel{cconstraint}
\frac{\eps {\gamma g_i n} +b_i}{r_i n-a_i} \le\eps\gamma\frac{ g_i}{r_i}+c_i.
\ee
We will appropriately define non-negative  $a_i$,   $b_i$ and $c_i$ with the following constraints\nek{:}
\bel{abconstraints}
 a_i< r_in/2, \qquad b_i\le \eps {\gamma g_i n}, \qquad c_i\le {\eps \gamma g_i /r_i \le 1/2,}
\ee
Note that requiring $c_i\le\eps \gamma g_i /r_i  \le 1/2$ ensures that we satisfy $f(i\eps)\le 1$.  So our definitions of these numbers just need to satisfy~\eqn{cconstraint} and~\eqn{abconstraints} for appropriate $\eps$, and allow (A1) and (A2) to hold with $i$ replaced by $i+1$. 
At this point we add the requirement that 
\bel{epsdef}
 \eps \sim n^{-\alpha},  \quad\mbox{where $0< \alpha<1/3$ is fixed}, 
\ee
with further conditions on $\eps$ to be imposed later, usually indirectly via conditions on   
$\alpha$.  Note that since $m=\Omega(n)$ and $\alpha<1/3$, for any such $\alpha$ we can always find such an $\eps$ for which $\eps m$ is an integer. One implication we will use is that 
\bel{epssquared}
\eps^2\ge 1/n
\ee
for $n$ sufficiently large.

In order to show that condition (A1) is satisfied after the $(i+1)$-st iteration, we need to estimate $|M(i+1)|$ for any $M\in \M^j(i)$, where $j > i+1$. First, we bound the number of edges  in $M$  that have at least one end vertex condemned  in step  (i) or (ii). We also call such edges {\em condemned}. Given $uv\in M(i)$, we know that the probability that $u$ (or $v$) is condemned after step (ii) is $f(i\eps)$. However, while the probability that a vertex is condemned in iteration $(i+1)$ is the same for all surviving vertices, vertices are not condemned independently. The following lemma shows that the probability that both $u$ and $v$ will be condemned is $O(f(i\eps)^2)$.  Note: the constants implicit in our $O(\cdot)$ notation are absolute. In the interest of continuity, we state the lemmas we need to prove \tref{thm:simple} below, and discuss some aspects of their proofs, but defer their proofs to \sref{subsec:lemproofs}.

\begin{lemma}~\lab{lem:del2} 
If $u$ and $v$ are distinct vertices in $G(i)$, then the probability that both $u$ and $v$ are condemned in iteration $i+1$ is $O(f(i\eps)^2)$. 
\end{lemma}

From this lemma, the probability that the edge $uv$ is condemned is $2f(i\eps)+O(f(i\eps)^2)$. By linearity, the expected number of condemned edges in any given surviving matching $M(i)\in \M^j(i)$ in the $(i+1)$-st iteration is 
$\big(2f(i\eps)+O(f(i\eps)^2)\big)|M(i)|$. Next, we address the effect of vertex collisions on the size of the  surviving  matchings. The following two lemmas bound the expected number of vertex collisions, and size of $\Phi(i+1)$, respectively. 

\begin{lemma}\lab{lem0:collision}
 Let $X_u$ be the number of edges incident with $u$ that are chosen in step {\rm (i)}, and let $Y_u=X_uI_{X_u\ge 2}$. 
With probability $1-o(\eps)$, 
\[
Y_u =  O\left(\max\left\{\frac{{\eps g_i } }{r_i^2n}d_u^{i+1}(i),\log^2 {n}\right\}\right).
\]
\end{lemma}

\begin{lemma}\lab{lem:collision}  With probability $1-o(\eps)$, we have
$|\Phi(i)|=O(\eps f(i\eps) {m} +\sqrt{\eps {m}}\log {n} )$. 
\end{lemma}

Thus, the treatment of vertex collisions does not change the size of each matching obtained from step (ii) significantly. The number of edges that are condemned but do survive, or are not condemned in steps (i) and (ii) but are deleted in step (iii), is bounded by $O(|\Phi(i)|)$.
It also follows from \lref{lem:collision} that step (iii) will not usually fail, as the number of matchings to be treated in that step is usually of much smaller order than $r_i n$, the approximate size of each matching.

Using such considerations, we are able to show that with high probability, the size of each surviving matching is concentrated around its expectation.

\begin{lemma}\lab{lem:A1}
With probability $1-o(\eps)$, for every $M\in \M^j(i)$ and $j > i+1$, we have
\[
|M(i+1)| = \big(1-2f(i\eps)+O(f(i\eps)^2)\big)|M(i)| +O( \eps f(i\eps) {m} + \sqrt{\eps {m}}\log {n}).
\]
\end{lemma} 
This provides us with enough information to specify $a_{i+1}$ as required for  (A1) after    iteration $i+1$.

Next we consider (A2). This requires us to  bound $d_v^j(i+1)$ for all $j > i+1$. Recall from \sref{sec:heuristic} that $E_v$ denotes the set of edges in $G$ that are incident with $v$. Let ${\mathcal E}=E_v\cap E(\M^j(i))$. Since every vertex in $G(i)$ is condemned  with probability $f(i\eps)$, by linearity, the expected number of  edges in ${\mathcal E}$ that are condemned is $f(i\eps) d_v^j(i)$, if $v$ survives after the $i$-th iteration. Again, \lref{lem:collision} ensures that the effect from vertex collision is small. This yields the following lemma.
\begin{lemma}\lab{lem:A2}
With probability $1-o(\eps)$, for every $v\in G(i+1)$ and $j > i+1$, 
\[
d_v^j(i+1) \le (1-f(i\eps))d_v^j(i)+O\left(\frac{\eps { g_i }}{r_i^2 n}\right)d_v^{i+1}(i)+O( \sqrt{\eps {m}}\log {n}).
\]
\end{lemma}
This lemma is strong enough for us to choose $b_{i+1}$ appropriately.

 We are now ready to complete the proof of \tref{thm:simple}.  We first write the requirements for (the $i+1$ versions of) (A1), (A2) and~\eqn{cconstraint} to be satisfied, using the inductive hypothesis, and then determine  $a_{i+1}$, $b_{i+1}$ and $c_i$
so as to satisfy these requirements as well as~\eqn{abconstraints}. We have
\bel{Meq}
|M(i+1)|- r_{i+1}n= T_1+T_2+T_3
\ee
where
\bean
T_1&= |M(i+1)|-|M(i)|= -2f(i\eps)|M(i)|+O\big(f(i\eps)^2|M(i)| + \eps f(i\eps) {m} + \sqrt{\eps {m}}\log {n}\big),\\
&\hspace{12cm} \mbox{(by \lref{lem:A1})}\\
T_2&= |M(i)|-r_in,\\
T_3&= (r_{i}-r_{i+1}) n =  2\eps\gamma n g_i  + O(\eps^2 \gamma^2n),
\end{align*}
where the last equation holds since  $r'(x)=-2\gamma g(x)$ and  $r''(x)=2\gamma^2$.
Now 
\[
 -2f(i\eps)|M(i)| +2\eps\gamma ng_i   = -2 f(i\eps)T_2 - 2c_ir_in
\]
and hence~\eqn{Meq} gives
\[
 |M(i+1)|- r_{i+1}n =   \big(1-2 f(i\eps)\big)T_2 - 2c_ir_in +O(f(i\eps)^2|M(i)|+ \eps f(i\eps) {m} + \sqrt{\eps {m}}\log {n}+\eps^2 n).
\]
It is an easy observation that (A1) implies $|T_2|\le a_i$, and this, along
with~\eqn{abconstraints} yields $f(i\eps)^2|M(i)|=O({\eps^2n g_i^2/r_i })=  O(\eps^2n) $. Indeed, $f(i\eps)= O(\eps { g_i/r_i })$ and thus, (A1) is satisfied after iteration $i+1$ provided that we define 
\bel{arequired}
a_{i+1}
= C_0\bigg(\eps^2{ g_i m/r_i } + \sqrt{\eps m}\log n  \bigg)
+  2c_ir_in+a_i\bigg(1-{\frac{2\eps\gamma g_i }{r_i }}\bigg),
\ee
where $C_0$ is a sufficiently large constant (subsequently to have a further condition imposed on it). Note that $\eps^2 n$ is absorbed by $\eps^2{ g_i m/r_i }$.

For (A2), we first rewrite
$$
d_v^j(i+1)-\eps {\gamma g_{i+1} n}  =  d_v^j(i+1)-\eps {\gamma g_i n}(1-f(i\eps))+\eps {\gamma n}\big(g_i-g_{i+1} - g_{i}  f(i\eps)\big).
$$
By the definition of $g$, we have ${g_i-g_{i+1}=\eps\gamma }$. Also, $ g_{i}  f(i\eps)\ge {\eps \gamma g_i^2 /r_i =\eps \gamma } $ as $c_i\ge 0$. Hence, using \lref{lem:A2} to bound the value of $d_v^j(i+1)$ occurring in the right hand side, we have (using $d_v^{i+1}(i)=O(\eps { g_i n})$ by (A2) and~\eqn{abconstraints})
\[
d_v^j(i+1)-\eps \gamma g_{i+1}n \le  (1-f(i\eps))\big(d_v^j(i)  -{\eps \gamma g_i n}\big) +O\big( {\eps^2 g_i^2 /r_i^2}+ \sqrt{\eps m}\log n \big).
\]
Thus, (A2) is satisfied after iteration $i+1$ provided that for a sufficiently large constant $C_0$,
\bel{brequired}
b_{i+1} =  (1-f(i\eps))b_i +C_0\big({\eps^2 g_i^2 /r_i^2}+ \sqrt{\eps m}\log n \big).
\ee
We choose $C_0$ sufficiently large to satisfy the bounds on it implied in deriving both \eref{arequired} and \eref{brequired}.
As part of the induction we are going to ensure the following strengthening of the constraints on $a_i$ and $b_i$  in~\eqn{abconstraints} (justified below):
\bel{asmall}
a_i\le \xi r_i n, \quad b_i \le \xi \eps { g_i n} \quad\mbox{for some fixed function $\xi = \xi(n)\to 0$}.
\ee
Then it follows that~\eqn{cconstraint} is satisfied for $n$ sufficiently large, provided we choose
\bel{crequired}
c_i= {\frac{\eps \gamma a_i g_i n  (1+2\xi)}{r_i^2n^2}}+\frac{b_i(1+2\xi)}{r_i n} \le  {\frac{\eps \gamma a_i g_i   (1+2\xi)}{r_i^2n}}+{\frac{2\xi \eps\gamma g_i }{r_i  }}. 
\ee

To complete the induction to the end of step $\tau-1$, it only remains to check the growth rates of $a_i$, $b_i$ and $c_i$ and see that they satisfy~\eqn{abconstraints} and~\eqn{asmall} (for an appropriate $\xi$), which can be assumed for smaller values of $i$ by induction.

Plugging~\eqn{crequired} into~\eqn{arequired}  and using \eqn{asmall} 
we get
\bel{rec}
 a_{i+1} \le A_i+B_ia_i
\ee
where 
\begin{align}
A_i &=  C_0\big( \eps^2{g_i m/r_i } + \sqrt{\eps m}\log n  \big) 
+  4b_i,\lab{e:Ai}\\
B_i &=   \frac{2\eps \gamma g_i   (1+2\xi)}{r_i} + 1 
-\frac{2\eps\gamma g_i }{r_i }  
= 1+  \frac{4\xi\eps\gamma g_i }{r_i }. \nonumber
\end{align}

If we turn the inequality~\eqn{rec}  into an  equality, we obtain a recurrence whose solution, from initial condition $a_0=0$, is  easily solved, and thus (since all coefficients are positive) implies
\bel{asol}
a_{i} \le \sum_{j=0}^{i-1} A_j\prod_{k=j+1}^{i-1} B_k.
\ee


Recall that the number of iterations the algorithm takes is $\tau=\lceil 1/\eps\rceil$. 
 For any $i\le\ceil{1/\eps}-1$,
\bea
\prod_{k=j+1}^{i-1} B_k &\le\exp\bigg((4\eps\gamma\xi )\sum_{k=j+1}^{i-1}\frac{1}{1-\eps \gamma  k} \bigg)\nonumber\\
&=\exp\bigg( (4\xi+o(1))
\int_{(\eps \gamma ) j}^{(\eps \gamma ) i  }  (1-x)^{-1} dx  \bigg) \nonumber\\
&\le{(1-\gamma)}^{ o(1)}=n^{o(1)} \lab{Bbound}
\end{align}
since  $\gamma=1-n^{-c}$.
 
We have by iterating~\eqn{brequired} (ignoring the negative term, which {turns out to give} no significant help) that
\bea 
b_i&\le b_0 + iC_0\big({\eps^2 g_i^2 /r_i^2}+ \sqrt{\eps m}\log n  \big)\nonumber\\
&\le ( 1 +iC_0)\sqrt{\eps m}\log n +iC_0{\eps^2/r_i},  \lab{bbound}
\end{align} 
 recalling that $b_0= \sqrt{\eps\gamma n}\log n$  and, as observed at the start of \sref{s:proof}, $m=\floor{\gamma n^{1+\delta}}$.
 This easily establishes the bound on $b_i$ in~\eqn{asmall} as long as
 \be
 \frac{1}{2}> \frac{3}{2}\alpha+\frac{\delta}{2}+c,\quad c<1/3. \lab{bcond}
 \ee
Now we turn to $a_i$.
Substituting~\eqn{bbound} into \eref{e:Ai} gives
\bel{Aj}
 A_j =O \big(\eps^2{ g_j m/r_j } + j\sqrt{\eps m}\log n+j{\eps^2 /r_j}  \big).
\ee
 Using this and~\eqn{Bbound} in~\eqn{asol}, and the bound 
 $i\le \ceil{1/\eps}-1$  gives
\bean
a_i& \le n^{o(1)}\cdot O\left(  \sum_{j=0}^{i-1} \eps^2{ g_j m/r_j }
+  \sum_{j=0}^{i-1}  j\sqrt{\eps m}\log n
+  \sum_{j=0}^{i-1}j{\eps^2 /r_j} \right).
\end{align*}
We can approximate $\sum_{j=0}^{i-1} \eps^2{ g_j m/r_j }$ and $\sum_{j=0}^{i-1}j{\eps^2 /r_j}$  as follows:
\[
\sum_{j=0}^{i-1} \eps^2{ g_j m/r_j } =O\left(\eps m\int_{0}^{(\eps\gamma)i} \frac{1}{1-x} dx\right)= O(\eps m\log(1/(1-\gamma)))=O(\eps m\log n),
\]
and
\[
 \sum_{j=0}^{i-1}j{\eps^2 /r_j}=O\left(\int_{0}^{(\eps\gamma )i} \frac{x}{(1-x)^2}dx\right) = O(1/(1-\gamma))=O(n^c),
\]
It then follows that
\begin{align*}
a_i 
& =n^{o(1)}\cdot O\left(\eps m +\sqrt{m}\eps^{-3/2}+n^c\right)
\end{align*}
as the logarithmic factors are absorbed by $n^{o(1)}$.
Since  $\tau=\ceil{1/\eps}$, we have  $\tau-1\le 1/\eps$,  and thus 
$r_{\tau-1}\ge (1-\gamma )^2=n^{-2c}$. As $r$ is monotonically decreasing, and  recalling that $\gamma=n^{-\delta}$,   $m\sim n^{1+\delta}$,  and $\eps\sim n^{-\alpha}$ from~\eqn{epsdef}, the above estimate for $a_i$ implies the bound for $a_i$ required in~\eqn{asmall}, provided that
$$
1-2c> \max\big\{-\alpha+1+\delta,(1+\delta+3\alpha)/2, c\big\}.
$$
As mentioned before, the first two bounds in~\eqn{abconstraints} follow from~\eqn{asmall}. The upper bound on $c_i$ in~\eqn{abconstraints} follows immediately from its definition in~\eqn{crequired}, in view of~\eqn{asmall}. 
Also, since $g_i\le 1$, we have the (final) upper bound, $1/2$, in~\eqn{abconstraints} provided 
$$
2c<\alpha+\delta.
$$
In summary, if these last two inequalities hold, as well as~\eqn{bcond}, then we have~\eqn{abconstraints} and~\eqn{asmall} . These three inequalities follow if we ensure that
\be
\delta+2c-\alpha<0, \quad \frac{3}{2}\alpha+\frac{\delta}{2}+2c<\frac{1}{2},\quad c<1/3. \lab{cond1}
\ee
By the theorem's hypothesis that $c<(1-4\delta)/10$,   there exists {$\alpha$} satisfying these conditions as well as the original $\alpha<1/3$ from~\eqn{epsdef}. (Note that the   bound   $c<1/3$   already follows from the theorem's hypothesis.)
We conclude that (A1) and (A2) are satisfied by induction, and hence with probability $1-o(1)$, the algorithm runs successfully to the end of the second-last iteration. Moreover, at the beginning of the last iteration,  each surviving matching has size at least $r_{\tau-1} n-a_{\tau-1} \ge r_{\tau-1}n/2$ by~\eqn{asmall}.

 Now we  argue that with probability $1-o(1)$, the algorithm finds a full rainbow matching in the last iteration. 
 The first inequality in~\eqn{cond1} gives  
\be 
2\eps m \le \frac{n^{1-2c}}{2} \le \frac{r_{\tau-1} n}{2}
\lab{cond3}
\ee
 for large $n$.
  There are at most $\eps m$ matchings remaining in the last iteration. So we can greedily choose one edge from each matching sequentially, since  $2\eps m \le r_{\tau-1}n/2$, by~\eqn{cond3}. \qed

\subsection{Proofs of lemmas}\label{subsec:lemproofs}
\textbf{Proof of Lemma \ref{lem:del2}.}
 %
Vertices $u$ and $v$ are both marked in step (i) if either
\bean
&\mbox{$uv\in E(\M^{i+1}(i))$ and $uv$ is chosen; or}\\
&\mbox{one edge in $E(\M^{i+1}(i)) \cap E_u$ is chosen and another edge in $E(\M^{i+1}(i)) \cap E_v$ is chosen}. 
\end{align*}
This probability is at most 
\be
\frac{1}{r_in -a_i}+\frac{d_u^{i+1}(i)}{r_in -a_i } \cdot \frac{d_v^{i+1}(i)}{r_in -a_i }=O \left(\eps\frac{g_i}{r_i}\right)^2=O(f(i\eps)^2), \lab{2vertices} 
\ee
because of~\eqn{abconstraints} and~\eqn{epssquared}, 
which imply that $1/(r_in-a_i)=O(\eps^2 g_i^2/r_i^2)$.
Vertices $u$ and $v$ are both condemned (marked or zapped) 
after step (ii) if and only if 
\bean
&\mbox{they are both marked in step (i); or}\\
&\mbox{one is condemned, and the other is zapped in step (ii)}. 
\end{align*}
We have shown the probability of the first case is $O(f(i\eps)^2)$. The probability of the second case is at most
\[
f(i\eps) P_{i+1}(v) + f(i\eps) P_{i+1} (u) = O(f(i\eps)^2).
\]
This is because {the probability of condemning $u$ is at most $f(i\eps)$ and conditional on $u$ being condemned and $v$ not being killed (with probability at most 1), the probability that $v$ is zapped is at most $P_{i+1}(v)$, as vertices are zapped independently in step (ii)}. The lemma follows.\qed \ss

\smallskip
\noindent
\textbf{Proof of Lemma \ref{lem0:collision}.} 
Recall that  $X_u$ denotes the number of edges incident with $u$ that are chosen in step (i), and $Y_u=X_uI_{X_u\ge 2}$. Immediately we have $X_u-1\le Y_u\le X_u$. 
Note that 
\[
Y_u\le X_u(X_u-1)= \sum_{x,y\in E_u\cap E(\M^{i+1}(i))} I_xI_y,
\]
where $I_x$ is the indicator variable that $x$ is chosen, and the summation is over all ordered pairs $(x,y)$. For each $u\in G(i)$, $u$ is incident with $d_u^{i+1}(i)\le {\eps \gamma g_i n+b_i}$ edges in $E(\M^{i+1}(i))$. Note that all edges in $E_u\cap E(\M^{i+1}(i))$ must belong to different matchings and therefore $\{I_x: x\in E_u\cap E(\M^{i+1}(i))\}$ are independent variables. Thus, the probability that any given $x$ and $y$ are both chosen is at most
\[
(r_in-a_i)^{-2}.
\]
Hence,
\[
\ex Y_u \le \big(d_u^{i+1}(i)\big)^2 (r_in-a_i)^{-2} = O\left(\frac{{\eps  g_i }}{r_i^2n}\right)d_u^{i+1}(i).
\]
It follows immediately that 
\[
\ex X_u\le 1+ O\left(\frac{{\eps  g_i }}{r_i^2n}\right)d_u^{i+1}(i),
\]
as $X_u\le 1+Y_u$.
Note that $X_u=\sum_{x\in E_u\cap E(\M^{i+1}(i))} I_x$, which is the sum of independent indicator variables. 
Applying Chernoff's bound, we obtain that with probability $1-o(\eps)$,
\[
Y_u\le X_u\le \max\{2\ex X_u, \log^2 {n}\}=O\left(\max\left\{\frac{{\eps g_i } }{r_i^2n}d_u^{i+1}(i),\log^2 {n}\right\}\right), \quad \forall u\in G(i).
\]
The lemma follows. \qed
\smallskip

\noindent\textbf{Proof of Lemma \ref{lem:collision}.} Let $Y=\sum_{u\in G(i)} Y_u$ where $Y_u$ is defined as in lemma~\ref{lem0:collision}. Then $|\Phi(i+1)|\le Y$. Thus it immediately follows that
\[
\ex |\Phi(i+1)|\le \ex Y =\sum_u \ex Y_u  = O\left(\frac{{\eps  g_i }}{r_i^2n}\right)\sum_u d_u^{i+1}(i).
\]
By (A1), $\sum_{u}d_u^{i+1}(i) \le 2(r_i n+a_i)\cdot \eps {m}$. Thus,
\[
\ex |\Phi(i+1)| = O\left(\frac{\eps^2 {m } g_i}{r_i}\right)=O(\eps f(i\eps) {m}). 
\]
Apply Azuma's inequality to $|\Phi(i+1)|$. Changing the choice $x$ to another edge $y$ in a matching $M$ would affect $|\Phi(i+1)|$ by at most $3$. To see this, let $x=uv$ and $y=u'v'$. A matching $M'$ that was in $\Phi(i+1)$ could be removed after changing $x$ to $y$, if $z\in M'$ was chosen, and $z$ is the only chosen edge, besides $x$, that is incident with $u$ (or $v$). There can only be at most two such matchings. So changing $x$ to $y$ would decrease $|\Phi(i+1)|$ by at most $3$, counting $M$ itself. Similarly, changing $x$ to $y$ would increase $|\Phi(i+1)|$ by at most $3$. Thus, by Azuma's inequality with Lipschitz constant 3,  we have that with probability $1-o(\eps)$, 
$|\Phi(i+1)|=\ex |\Phi(i+1)| + O(\sqrt{\eps {m}}\log {n}) = O(\eps f(i\eps) {m}+\sqrt{\eps {m}}\log {n})$.
 \qed \medskip

We will use the following Azuma-Hoeffding inequality to prove concentration of various variables.
\begin{thm}[\cite{Azuma,Hoeffding}]\lab{thm:Azuma}
Let $X_0,X_1,\ldots$ be a martingale satisfying $|X_i-X_{i-1}|\le \delta_i$ for every $i\ge 1$. Then, for every $t\ge 0$,
\[
\pr\big(|X_n-X_0|\ge t\big)\le 2\exp\left(-t^2\Big/2\sum_{i=1}^n \delta_i^2\right).
\]
\end{thm}

\noindent\textbf{Proof of Lemma \ref{lem:A1}.} 
We have argued that 
\be
\ex\big(|M(i+1)|\,\big|\, G(i)\big) = \big(1-2f(i\eps)+O(f(i\eps)^2)\big)|M(i)|+O(\eps f(i\eps) {m} +\sqrt{\eps {m}}\log {n}),\lab{eq:M}
\ee
where the main term comes from considering the edges that are condemned, and the error term accounts for a correction term due to vertex collision,  by \lref{lem:collision}.

For concentration, first consider $X$, the number of edges chosen in $\Psi(i+1)$ in $M\in \M^j(i)$ in step (i). Let $\ex X=Y_1,\ldots, Y_{\eps m}=X$ be the Doob's martingale constructed by the conditional expectation of $X$ under the edge exposure process where edges in $\Psi(i+1)$ are revealed sequentially. Apply \tref{thm:Azuma} to the martingale $(Y_i)$. It is easy to see that changing a single edge $x\in \Psi(i+1)$ to another edge $y$ would change $X$ by at most two. Thus, the probability that $X$ deviates from $\ex X$ by more than $t=\sqrt{\eps m} \log n$ is at most
$2\exp(-t^2/8\eps m)=o(n^{-2})$. Taking the union bound over all $M\in \M^j(i)$ we obtain the desired deviation in the lemma with probability at least $1-o(\eps)$. Next, we consider the number of edges zapped in $M$ in step (ii). Condition on the set of edges that survive step (i). Each surviving vertex $u$ is zapped independently with probability $P_{i+1}(v)$. For each $M\in \M^j(i)$ consisting of edges surviving after step (i), the $2|M|$ vertices incident with $M$ are independently zapped with probabilities all bounded by $f(i\eps)$. Let $Y$ denote the number of vertices zapped. Then, $Y$ is the sum of at most $2n$ independent Bernoulli variables. By the Chernoff-Hoeffding bound~\cite[Theorem 1.1]{dp}, the probability that $Y$ deviates from its expectation by more than $\sqrt{f(i\eps) n} \log n$ is at most $n^{-2}$. Taking the union bound over all $M$, again with probability at least $1-o(\eps)$ we have the desired deviation as in the lemma. 
Finally, the change to $|M(i)|$ due to step (iii) is absorbed by the error term in~\eqn{eq:M} by \lref{lem:collision}.\qed\ss

\noindent\textbf{Proof of \lref{lem:A2}.} With arguments similar to the proof for \lref{lem:A1}, we can apply \tref{thm:Azuma} to prove concentration for the number of neighbours of $v$ (in chunk $j$) condemned in step (i) and then for the number of neighbours zapped in step (ii) using  the Chernoff-Hoeffding bound. 
By \lref{lem0:collision}, vertex collision will affect $d_v^j(i+1)$ by $O\left(\eps {\gamma g_i /r_i^2 n}\right)d_v^{i+1}(i)$. The treatment of vertex collisions in step (iii) can only decrease (the bound on)  $d_v^j(i+1)$. The lemma follows.\qed

\section{Multigraphs}
 
The proof for \tref{thm:multiple} follows almost exactly that of \tref{thm:simple} with $\delta=0$. We run the same randomised algorithm with the same parameters $g_i$ and $r_i$, but with different $a_i$, $b_i$ and $c_i$. The reason is that due to the multiplicities of the multiple edges, variables $|M(i)|$ and $d_v^j(i)$ are not as concentrated as in the simple graph case and thus we expect larger $a_i$, $b_i$ and $c_i$ here. We briefly sketch the proof.
Now we assume that $m=\floor{\gamma n}$ where $\gamma=1-\eps_0$ and $\eps_0>0$ is an arbitrarily small constant, and $\eps>0$ is going to be a constant that depends on $\eps_0$. Let $\mu$ denote the maximum multiplicity of the multiple edges in $G$. Note also that here $m=\Theta(n)$.

\lref{lem:base} holds in the multigraph case with the same $b_0$.
For \lref{lem:del2}, the probability that an edge between $u$ and $v$ is chosen is bounded by $\mu/(r_i n-a_i)$. Thus, for~\eqn{2vertices} to hold, we require
 \[
 \frac{\mu}{r_in-a_i}
=O\left(\frac{\eps^2(g_i m)^2}{(r_i n)^2}\right),
 \]
 which is guaranteed if we assume
\be
\eps^2\ge \mu /n.\lab{del2New}
\ee
 Thus, \lref{lem:del2} holds after replacing the condition $\eps^2\ge 1/n$ by~\eqn{del2New}.
 
 Lemmas~\ref{lem0:collision},~\ref{lem:collision} and~\ref{lem:A1} hold as they are. 
 For \lref{lem:A2}, note that deleting a single vertex (both in steps (i) and (ii)) can alter $d_v^j(i+1)$ by $\mu$. Therefore, the Lipschitz constant becomes $\mu$. Applying Azuma's inequality to both steps (i) and (ii), we deduce \lref{lem:A2} where $O(\sqrt{\eps m}\log n)$ is replaced by $O(n/\log n)$, if $\mu=O(\sqrt{n}/\log^2 n)$.  
 
 These lead to recursions for $a_i$ as in~\eqn{rec}, and for $b_i$ as 
   \[
  b_{i+1} = (1-f(i\eps)) b_i+O(n/\log n).
  \]
  Immediately we have $b_i=O(in/\log n)$.
  Substituting into~\eqn{rec}, we have
  \[
  a_{i+1}\le \left(1+\frac{4\xi\eps\gamma g_i }{r_i }\right)a_i + O(\eps^2 n g_i/r_i +in/\log n).
  \]
  Solving the recursion as before we get
  \bean
  a_i&=\eps_0^{o(1)}\cdot O\left(\sum_{j=0}^{i-1} \eps^2 n g_j/r_j +\sum_{j=0}^{i-1}jn/\log n \right)\\
  &=\eps_0^{o(1)}\cdot O\left(\eps n + n/\eps\log n\right).
  \end{align*}
  Hence there exists a constant $C>0$ such that $a_i\le C(\eps n+n/\eps \log n)$ for all $1\le i\le 1/\eps$.
  Let  $i_1=\lceil 1/\eps\rceil -1$. Then $r_{i_1}\ge (1-\gamma)^2=\eps_0^2$. Choose $\eps>0$ sufficiently small such that 
  \[
  C_0(\eps n+n/\eps \log n) \le \frac{\eps_0^2}{4}n.
  \]
  Then, with the same argument as before, a greedy search in the last iteration of the randomised algorithm succeeds in finding a full rainbow matching in $\M$
with high probability.
  
  \section{Hypergraphs}
  
  The proof of \tref{thm:hypergraph} is again similar.  Let $G$ be a $k$-uniform hypergraph. We give a quick sketch here and just point out the differences. The randomised algorithm extends to hypergraphs in a natural way. Thus, every vertex is deleted in the $(i+1)$-st iteration with probability 
  \[
  f(i\eps)\approx \frac{\eps  \gamma  g(i\eps)}{r(i\eps)},
  \]
  where $\gamma=1-\eps_0$.
Now every hyperedge in a matching is deleted with probability approximately
  $kf(i\eps)$, as there are $k$ vertices in a hyperedge, For each surviving vertex $v$, each incident hyperedge is deleted with probability approximately $(k-1)f(i\eps)$, as this hyperedge is deleted if one of the other $k-1$ vertices contained in it is deleted. Hence, we find that $r(x)$ and $g(x)$ obey the following differential equations
  \[
  r'=-k \gamma g(x), \ \ g'(x)= - (k-1)\gamma \frac{g(x)^2}{r(x)},
  \]
  with initial conditions $r(0)=1$ and $g(0)=1$. The solution to these differential equations is
  \[
  r(x)=(1-\gamma x)^k,\ \ g(x)=(1-\gamma x)^{k-1}.
  \]
  The proof that $|M(i)|$ and $d_v^j(i)$ are concentrated around $r_i n$ and $\eps \gamma g_i n$ follows in the same manner as in \tref{thm:simple}. Lemmas~\ref{lem:del2} and~\ref{lem:A2} need to be modified as in \tref{thm:multiple}. Here, the affect of codegrees, i.e.\ the maximum number of hyperedges containing a pair of vertices plays the same role of affecting the Lipschitz constants as the maximum multiplicity in \tref{thm:multiple}. This yields \tref{thm:hypergraph}.\qed

\section{Counterexamples to some conjectures on rainbow matchings}\label{s:CEs}

In this section we describe counterexamples to Conjectures 2.5 and 2.9
in \cite{ahbe},  as well as Conjectures 5.3, 5.4, 6.1 and 6.2 in
\cite{ABCHS}.  Before doing so, we need to describe how rainbow
matchings in graphs can be viewed as matchings in 3-uniform
hypergraphs.  Suppose that $G$ is a (not necessarily properly)
edge-coloured graph in which we are interested in finding a rainbow
matching. We make a 3-uniform hypergraph $H$ from $G$ as follows.
The vertices of $H$ are $V(G)\cup V_1$ where $V_1$ is the set of colours
used on edges of $G$. For each edge $\{u,v\}$ of $G$ with colour $c\in V_1$
there is a hyperedge $\{u,v,c\}$ in $H$.
Now a full rainbow matching in $G$ corresponds to a matching of
$H$ that covers all of the vertices in $V_1$.  If $G$ happens to be
bipartite with bipartition $V_2\cup V_3$,
then $H$ will be {\em tripartite}, because its vertices can
be partitioned as $V_1\cup V_2\cup V_3$ such that every hyperedge
includes one vertex from each of these three sets. 

\def\GG{\mathcal{G}}
\def\HH{\mathcal{H}}

\begin{figure}
  \[
  \includegraphics[scale=0.5]{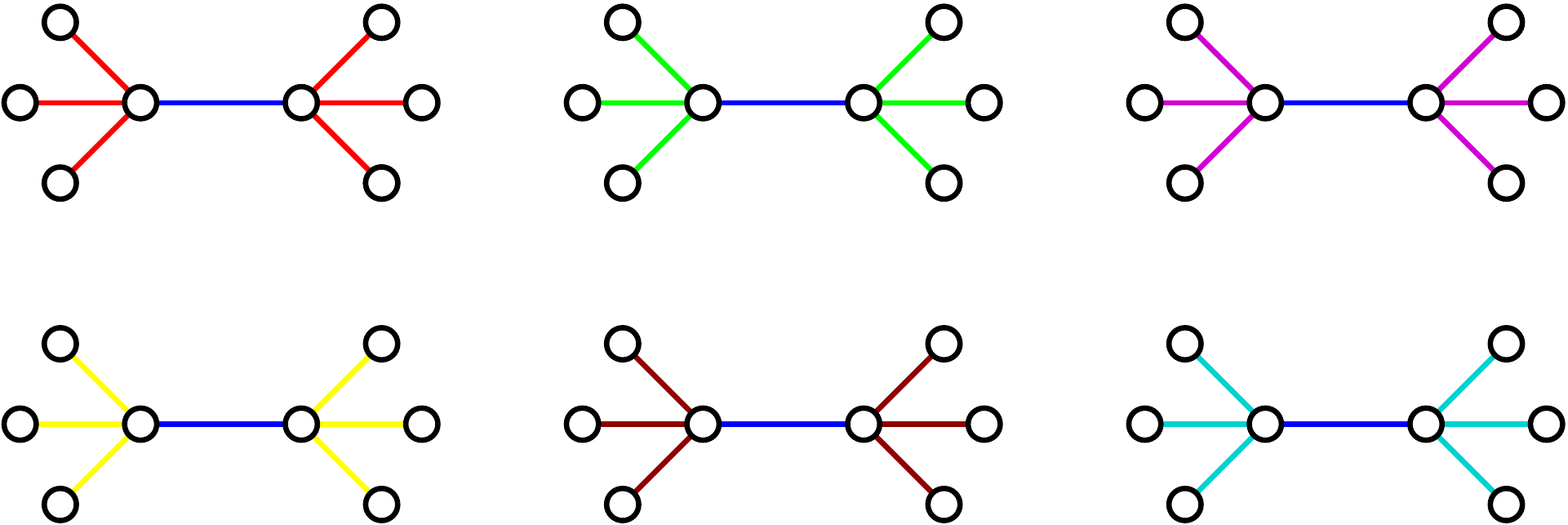}  
  \]
  \caption{\label{f:dblstar}The graph $\GG_6$}
\end{figure}

Let $m$ be a positive even integer.
We now construct a bipartite graph $\GG_m$ whose edges are (not properly)
coloured using $m$ colours in such a way
that there is no full rainbow matching. There are $m$ components in $\GG_m$,
each isomorphic to a double star which has two adjacent central
vertices each of which has $m/2$ leaves attached to it. The edge
between the central vertices in each double star is coloured blue. In
each component, the edges connected to leaves all have one
colour (not blue), which is specific to that component. Hence there are
$m+1$ colours overall, and each colour appears on $m$ edges.
\fref{f:dblstar} shows $\GG_6$.
There is
no full rainbow matching in $\GG_m$ because such a matching must include a blue
edge from some double star $S$. However, the colour of the other edges
in $S$ then cannot be represented in the matching.

Let $\HH_m$ be the tripartite hypergraph corresponding to $\GG_m$.
Let $V_1$ be the vertices of
$\HH$ corresponding to the colours, and $V_2,V_3$ the sets of vertices
corresponding to a bipartition of $\GG$.  Then every vertex in $V_1$ has
degree $m$. The vertices in $V_2\cup V_3$ all have degree either $1$
or $m/2+1$. Thus the (minimum) degree $\delta(V_1)$ of a vertex in
$V_1$ is nearly double the maximum degree $\Delta(V_2\cup V_3)$ of the
vertices outside $V_1$.  Interestingly, Aharoni and Berger
\cite[Thm~2.6]{ahbe} showed that in any tripartite hypergraph if
$\delta(V_1)\ge2\Delta(V_2\cup V_3)$ then there must be a
$|V_1|$-matching. Our hypergraph $\HH_m$ shows that their theorem is close
to tight. However, they made the following conjecture
\cite[Conj.~2.5]{ahbe} (repeated as \cite[Conj.~5.3]{ABCHS}, and rephrased
as \cjref{cj:falseconjbip} in our introduction).

\begin{con}\label{cj:falseconj}
  Let $H$ be a hypergraph with a vertex tripartition $V(H)=V_1\cup V_2\cup V_3$
  such that every hyperedge includes exactly one vertex from $V_i$ for $i=1,2,3$.
  If $\delta(V_1)>\Delta(V_2\cup V_3)$ then $H$ has a $|V_1|$-matching.
\end{con}

\begin{figure}
  \[
  \includegraphics[scale=0.4]{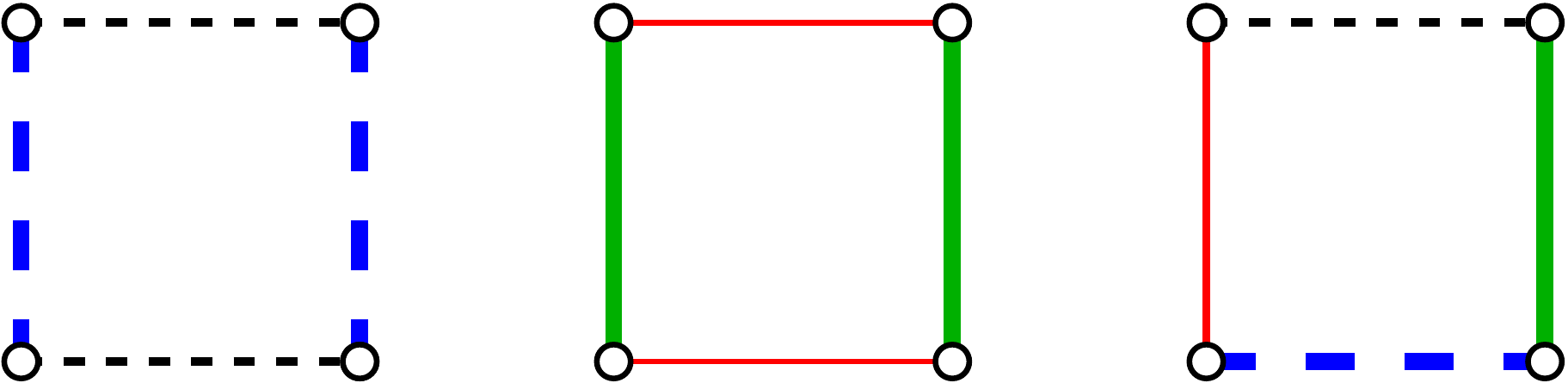}  
  \]
  \caption{\label{f:4cycleCE}A $2$-regular graph with no rainbow matching.}
\end{figure}

Note that $\HH_m$ disproves \cjref{cj:falseconj} whenever $m\ge4$.
Another counterexample to \cjref{cj:falseconj}
is based on the graph in \fref{f:4cycleCE}, which has no rainbow matching.
The corresponding tripartite hypergraph has 
$\delta(V_1)=3>2=\Delta(V_2\cup V_3)$. The line graph of the graph
in \fref{f:4cycleCE} was published in \cite{Alon92} and its complement
was published in \cite{Jin92}. In both cases the focus of the investigation
was slightly different from ours, so the generalisations that were offered
are not relevant for us.

Conjecture 2.9 of \cite{ahbe} generalises \cjref{cj:falseconj}, so it too is false.
Similarly, \cite[Conj.\,6.1]{ABCHS} asserts that
if $\delta(V_1)\ge2+\Delta(V_2\cup V_3)$ then there must be a
$|V_1|$-matching, so $\HH_m$ is a counterexample whenever $m\ge6$.

Finally, we consider Conjectures 5.4 and 6.2 from
\cite{ABCHS}. These deal with the case when the initial graph is not
necessarily bipartite, so the resulting hypergraph is not necessarily
tripartite. Nevertheless they consider full rainbow matchings
in an edge-coloured graph. Or equivalently, $|V_1|$-matchings in
a $3$-uniform hypergraph $H$ in which every hyperedge includes exactly one
vertex in the set $V_1$. The conjectures assert that such a matching
will exist provided that $\delta(V_1)\ge2+\Delta(V(H)\setminus
V_1)$. Again, $\HH_m$ provides a counterexample. Indeed, it
shows that the $2$ cannot be replaced by any constant.

  \let\oldthebibliography=\thebibliography
  \let\endoldthebibliography=\endthebibliography
  \renewenvironment{thebibliography}[1]{%
    \begin{oldthebibliography}{#1}%
      \setlength{\parskip}{0.4ex plus 0.1ex minus 0.1ex}%
      \setlength{\itemsep}{0.4ex plus 0.1ex minus 0.1ex}%
  }%
  {%
    \end{oldthebibliography}%
  }


\begin{thebibliography}{99}

\bibitem{ahbe} R.~Aharoni and E.~Berger, Rainbow matchings in $r$-partite $r$-graphs, \emph{Electron. J. Combin.} \textbf{16(1)} (2009) \#R119.

\bibitem{ABCHS}
R. Aharoni, E. Berger, M. Chudnovsky, D. Howard and P. Seymour,
Large rainbow matchings in general graphs,
arXiv:1611.03648v1. 

\bibitem{ahchho} R.~Aharoni, P.~Charbit and D.~Howard, On a generalization of the Ryser-Brualdi-Stein conjecture, \emph{J. Graph Theory} \textbf{78} (2015), 143--156.
  
\bibitem{ahkozi} R.~Aharoni, D.~Kotlar and R.~Ziv, 
Representation of large matchings in bipartite graphs,
{\it SIAM J.\ Discrete Math.} {\bf31} (2017), 1726--1731.

\bibitem{Alon92}
N.~Alon, 
The strong chromatic number of a graph,
\emph{Random Structures Algorithms} {\bf3} (1992), 1--7. 

\bibitem{Azuma} K. Azuma, Weighted sums of certain dependent random variables, {\em Tohoku Mathematical Journal}, Second Series 19.3 (1967): 357--367.

\bibitem{BGS}
J.~Bar\'at, A.~Gy\'arf\'as and G.\,N.~S\'ark\"ozy,
Rainbow matchings in bipartite multigraphs,
{\it Period. Math. Hungar.} {\bf74} (2017), 108--111. 
  
\bibitem{BW14}
J.~Bar\'at and I.\,M.~Wanless,
Rainbow matchings and transversals,
{\it Australas.\ J.\ Combin.\/} {\bf59}, (2014) 211--217.

\bibitem{BHWWW17}
D.~Best, K.~Hendrey, I.\,M.~Wanless, T.\,E.~Wilson and D.\,R.~Wood,
Transversals in Latin arrays with many distinct symbols,
{\em J.\ Combin.\ Des.} {\bf26} (2018), 84--96.


\bibitem{CW17}
N.\,J.~Cavenagh and I.\,M.~Wanless, 
Latin squares with no transversals,
{\it Electron.\ J.\ Combin.\/} {\bf24(2)} (2017), \#P2.45.

\bibitem{cleh} D.~Clemens and J.~Ehrenm\"{u}ller, An improved bound on the sizes of matchings guaranteeing a rainbow matching, \emph{Electron. J. Combin.} \textbf{23(2)} (2016), \#P2.11.

\bibitem{dp} D. P. Dubhashi, and A. Panconesi, Concentration of measure for the analysis of randomized algorithms, {\em Cambridge University Press}, 2009.

\bibitem{hajo} R.~H\"aggkvist and A.~Johansson, Orthogonal Latin rectangles, \emph{Combin. Probab. Comput.} \textbf{17} (2008), 519--536.

\bibitem{HS08} P.~Hatami and P.~W.~Shor, A lower bound for the length of a partial transversal in a Latin square, \emph{J. Combin. Theory Ser. A.} \textbf{115} (2008), 1103--1113.
 
\bibitem{Hoeffding} 
W. Hoeffding, Probability inequalities for sums of bounded random variables, 
{\em J. Amer. Statist. Assoc.} 58.301 (1963): 13--30.
 
\bibitem{Jin92}
G.~Jin, Complete subgraphs of $r$-partite graphs,
\emph{Combin. Probab. Comput.} {\bf1} (1992), 241--250.

\bibitem{KY17}
P.~Keevash and L.~Yepremyan,
Rainbow matchings in properly-coloured multigraphs,
arXiv:1710.03041.

\bibitem{KY18}
P.~Keevash and L.~Yepremyan,
On the number of symbols that forces a transversal,
arxiv:1805.10911.

\bibitem{kozi} D.~Kotlar and R.~Ziv, Large matchings in bipartite graphs have a rainbow matching, \emph{European J. Combin.} \textbf{38} (2014), 97--101.

\bibitem{McDiarmid} C. McDiarmid, Concentration for independent permutations, {\em Combin. Probab. Comput.} 11.2 (2002): 163--178.

\bibitem{MPS18} 
R. Montgomery, A. Pokrovskiy and B. Sudakov,
Decompositions into spanning rainbow structures,
arxiv:1805.07564.

\bibitem{Pok16} 
A. Pokrovskiy,
An approximate version of a conjecture of Aharoni and Berger,
arXiv:1609.06346.

\bibitem{Pok17} 
A. Pokrovskiy,
Rainbow matchings and rainbow connectedness,
{\it Electron.\ J.\ Combin.\/} {\bf24(1)} (2017), \#P1.13.

\bibitem{ry} H.~Ryser, 
Neuere probleme der kombinatorik, 
\emph{Vortr\"{a}ge \"{u}ber Kombinatorik, Oberwolfach} (1967), 69--91.
   
\bibitem{Ste75} S.K.~Stein, 
Transversals of Latin squares and their generalizations, 
\emph{Pacific J.~Math.} \textbf{59} (1975), 567--575.

\bibitem{transurv} I.\,M.~Wanless,
``Transversals in Latin squares: A survey'', in
R.~Chapman (ed.), {\it Surveys in Combinatorics 2011},
London Math. Soc. Lecture Note Series {\bf392}, 
Cambridge University Press, 2011, pp403--437.

\end{thebibliography}
\end{document}